\documentclass[article]{elsarticle}

\usepackage{lineno,hyperref}
\usepackage{amsmath}
\usepackage{amssymb}
\usepackage{graphicx}
\usepackage{epic}
\usepackage{bxeepic}
\usepackage{color}
\usepackage{xcolor}
\usepackage{stmaryrd}
\modulolinenumbers[5]

\begin{document}

\begin{frontmatter}

\title{Geometric multigrid methods for Darcy--Forchheimer flow in fractured porous media}

\author[a]{A. Arrar\'as}
\ead{andres.arraras@unavarra.es}
\author[b,c]{F.J. Gaspar}
\ead{fjgaspar@unizar.es}
\author[a]{L. Portero\corref{cor1}}
\ead{laura.portero@unavarra.es}
\author[b]{C. Rodrigo}
\ead{carmenr@unizar.es}
\cortext[cor1]{Corresponding author.}
\address[a]{Departamento de Estad\'{i}stica, Inform\'atica y Matem\'aticas, Universidad P\'ublica de Navarra, Edificio de Las Encinas, Campus de Arrosad\'{\i}a, 31006 Pamplona, Spain}
\address[b]{IUMA, Departamento de Matem\'atica Aplicada, Universidad de Zaragoza, Pedro Cerbuna 12, 50009 Zaragoza, Spain}
\address[c]{Centrum Wiskunde \& Informatica (CWI), 1098 XG Amsterdam, The Netherlands}

%
%
%

\begin{abstract}
In this paper, we present a monolithic multigrid method for the efficient solution of flow problems in fractured porous media. Specifically, we consider a mixed-dimensional model which couples Darcy flow in the porous matrix with Forchheimer flow within the fractures. A suitable finite volume discretization permits to reduce the coupled problem to a system of nonlinear equations with a saddle point structure. In order to solve this system, we propose a full approximation scheme (FAS) multigrid solver that appropriately deals with the mixed-dimensional nature of the problem by using mixed-dimensional smoothing and inter-grid transfer operators. Remarkably, the nonlinearity is localized in the fractures, and no coupling between the porous matrix and the fracture unknowns is needed in the smoothing procedure. Numerical experiments show that the proposed multigrid method is robust with respect to the fracture permeability, the Forchheimer coefficient and the mesh size.
\end{abstract}


\end{frontmatter}


\section{Introduction}\label{sec:intro}

Modeling and simulation of fluid flow in fractured porous media is a challenging task which is getting increasing attention in recent years, due to the wide range of applications in which plays an essential role. Different fracture models have been proposed in the last decades based on the spatial scale under consideration and the knowledge of the fracture distribution. 
On the one hand, double-continuum models are suitable for regularly distributed micro-fractures showing interconnections with the surrounding matrix. Such models assume the existence of a mass transfer function between the bulk and the fractures \cite{bar:zhe:koc:60}, and are usually derived via homogenization theory \cite{arb:dou:hor:90}. On the other hand, discrete fracture networks consider sets of individual macro-fractures which are isolated from the porous matrix \cite{for:fum:sco:ruf:14,pic:erh:dre:12}. Typically, these networks are obtained stochastically and provide information about the orientation, density, size and hydrological properties of the fractures \cite{ben:ber:pie:sci:14}. In these latter models, fluid exchange between the fractures and the matrix is not allowed, so that flow is restricted to the fracture network.
If we properly combine the preceding models, we may construct what we refer to as discrete fracture-matrix models: sets of individual macro-fractures, similar to those arising in discrete fracture networks, but suitably coupled with the surrounding matrix, as in double-continuum models. These are the models considered in this paper. More precisely, we suppose that fractures can be represented as $(n-1)$-dimensional interfaces immersed into an $n$-dimensional porous matrix, thus giving rise to the so-called mixed-dimensional or interface models \cite{kei:fum:ber:ste:17,Martin-Jaffre-Roberts}.

Most earlier works using this approach suppose that the flow within the fractures and in the porous matrix is described by Darcy's law \cite{fle:fum:sco:16,delpra:fum:sco:17,fri:mar:rob:saa:12}. Darcy's law has been shown to govern single-phase incompressible flow in porous media at specific flow regimes where the velocity is low. This is the case, for example, of subsurface reservoirs and aquifers, where a low permeability of the porous matrix implies low velocities. However, in proximity to wellbores or within high-permeability fractures, velocities are higher, thus requiring the use of alternative nonlinear flow models \cite{bal:mik:whe:10}. The simplest of such models is based on the addition of a quadratic correction term in the velocity to the linear Darcy model. The new model, referred to as Forchheimer's law, combines the contribution of viscous and inertial effects: at low flow rates, the viscous effect is dominant and the model reduces to Darcy's law; at increasing flow rates, however, the inertial effect gains relevance and plays a significant role \cite{gee:74}. Remarkably, other nonlinear correction terms --e.g., cubic \cite{mei:aur:91}, polynomial \cite{bal:mik:whe:10} or exponential \cite{pan:fou:06}-- have also been proposed in the literature.

The validity of Forchheimer's law in a certain range of velocities for laminar flow has been established empirically (see \cite{hua:ayo:08,bar:con:04} and references therein). From a theoretical viewpoint, the Forchheimer model has been deduced using homogenization methods \cite{che:lyo:qin:01,gio:97}, volume averaging \cite{rut:ma:92,whi:96}, and related techniques \cite{has:gra:87}. Existence, uniqueness and regularity results have been derived in \cite{ami:91,fab:89,kna:sum:16}. Numerically, different strategies --ranging from mixed finite elements \cite{dou:pae:gio:93,kie:16,par:05} to block-centered finite differences \cite{rui:pan:12,rui:pan:17} and multipoint flux approximation methods \cite{xu:lia:rui:17}-- have been applied to obtain approximate solutions of this model.

In this work, we are concerned with the numerical solution of a discrete fracture-matrix model which couples Darcy flow in the porous matrix with Forchheimer flow within the fractures. The solvability of this problem is analyzed in \cite{kna:rob:14}. In \cite{fri:rob:saa:06,fri:rob:saa:08}, numerical approximations are obtained using the lowest order Raviart--Thomas mixed finite elements in combination with a domain decomposition technique. In both works, the nonlinear system stemming from the Forchheimer equation is solved using fixed-point iteration and quasi-Newton methods. Alternative efficient solvers for various discretizations of the isolated Forchheimer model include the Peaceman--Rachford iteration scheme \cite{gir:whe:08}, different variants of the two-grid method \cite{rui:liu:15,sun:rui:18}, and a multigrid method based on the so-called full approximation scheme (FAS) \cite{hua:che:rui:18}. In the spirit of this latter work, we propose a monolithic mixed-dimensional multigrid method that extends our earlier work \cite{arr:gas:por:rod:18} for the Darcy--Darcy coupling to the Darcy--Forchheimer case. Note that, in this case, the mixed-dimensional approach establishes a connection between dimensionality and nonlinearity: an $n$-dimensional linear Darcy problem is coupled with an $(n-1)$-dimensional nonlinear Forchheimer problem. For the discretization, we consider a finite volume method that combines control volumes of different dimensions in the fractures and the porous matrix. The nonlinear system stemming from the discretization has a saddle point structure, and can be suitably handled using the FAS multigrid solver \cite{Bra77}.

Multigrid methods are well known to be among the fastest solvers for the solution of linear and nonlinear systems of equations, showing very often optimal computational cost and convergence behavior \cite{TOS01}. 
The performance of multigrid algorithms strongly depends on the choice of their components, so that many details are open for discussion and decision in the design of a multigrid method for a target problem. In the framework considered here, where a mixed-dimensional problem needs to be solved, it seems natural to combine two-dimensional smoothing and inter-grid transfer operators for the unknowns in the porous matrix with their one-dimensional counterparts within the fracture network. Regarding the smoother, due to the saddle point character of the resulting system, a Vanka-type relaxation is proposed for both the unknowns in the porous matrix and within the fractures. This class of smoothers was firstly proposed by Vanka in \cite{vanka} for the multigrid solution of the staggered finite difference discretization of the Navier--Stokes equations and, since then, it has been applied to different problems in both computational fluid and solid mechanics. In particular, here we consider  a standard two-dimensional five-point Vanka smoother for the unknowns in the porous matrix, and its one-dimensional three-point nonlinear counterpart for those unknowns within the fractures. The inter-grid transfer operators that act on the different unknowns are dictated by the corresponding one- and two-dimensional staggered location of the grid-points within the fractures and the porous matrix, respectively. The proposed mixed-dimensional multigrid method is shown to be robust with respect to the fracture permeability, the mesh size, and the so-called Forchheimer coefficient, which represents a measure of the strength of the nonlinearity (see problem \eqref{cont:problem}).

The rest of the paper is organized as follows. In Section \ref{sec:model}, we describe the discrete fracture-matrix model coupling Darcy flow in the porous matrix with Forchheimer flow in the fractures. The finite volume spatial discretization is formulated in Section \ref{sec:discretization}, where we further specify the resulting nonlinear system of algebraic equations. In Section \ref{sec:multigrid}, we introduce a monolithic mixed-dimensional multigrid method for solving such a system. Finally, we report a collection of numerical experiments in Section \ref{sec:experiments}, illustrating the robustness of the proposed method with respect to different parameters.

\section{The continuous problem}\label{sec:model}

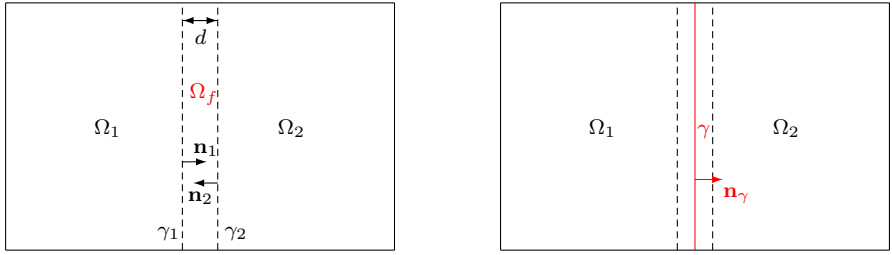
\begin{figure}[t]
	\begin{center}
		\unitlength=0.047cm
		\begin{picture}(260,100)
		\put(5,20){\line(1,0){110}}\put(115,20){\line(0,1){70}}\put(5,90){\line(1,0){110}}\put(5,20){\line(0,1){70}}
		\dashline[20]{2}(55,20)(55,90)
		\dashline[20]{2}(65,20)(65,90)
		\put(30,53){{\footnotesize$\Omega_1$}}
		\put(57,63){{\color{red}\footnotesize$\Omega_f$}} 
		\put(82,53){{\footnotesize$\Omega_2$}}
		\put(55,85){\vector(1,0){10}}
		\put(65,85){\vector(-1,0){10}}
		\put(58.5,78.5){{\footnotesize$d$}}
		\put(48,24){{\footnotesize$\gamma_1$}}
		\put(67,24){{\footnotesize$\gamma_2$}}
		\put(55,45){\vector(1,0){7}}
		\put(58,47.5){\footnotesize$\mathbf{n}_1$}
		\put(65,39){\vector(-1,0){7}}
		\put(56.5,34){\footnotesize$\mathbf{n}_2$}
		
		\put(145,20){\line(1,0){110}}\put(255,20){\line(0,1){70}}\put(145,90){\line(1,0){110}}\put(145,20){\line(0,1){70}}
		\put(200,20){\color{red}\line(0,1){70}}
		\dashline[20]{2}(195,20)(195,90)
		\dashline[20]{2}(205,20)(205,90)
		\put(170,53){{\footnotesize$\Omega_1$}}
		\put(200.7,53){{\color{red}\footnotesize$\gamma$}} 
		\put(222,53){{\footnotesize$\Omega_2$}}
		\put(200,40){\color{red}\vector(1,0){8}}
		\put(208,35){\color{red}\footnotesize$\mathbf{n}_{\gamma}$}
		\end{picture}
	\end{center}\vspace*{-1cm}
	\caption{Schematic representation of the original domain (left) and the reduced domain (right).}
	\label{fig:one:fracture}
\end{figure}

Let $\Omega\subset\mathbb{R}^2$ be an open, bounded, and convex polygonal domain, whose boundary is denoted by $\Gamma=\partial\Omega$. We consider a single-phase incompressible flow in $\Omega$ governed by the mass conservation equation, together with Forchheimer's law that relates the gradient of the pressure $p$ to the flow velocity $\mathbf{u}$, i.e.,
\begin{equation}\label{cont:problem}
\begin{alignedat}{2}
\left(1+\beta\, |\mathbf{u}|\right)\mathbf{u}&=-\mathbf{K}\nabla p\hspace*{1.5cm}&&\mbox{in }\Omega,\\
\nabla\cdot\mathbf{u}&=q\hspace*{1.5cm}&&\mbox{in }\Omega,\\
p&=0\hspace*{1cm}&&\mbox{on }\Gamma.
\end{alignedat}
\end{equation}
Here, $\beta$ represents the dynamic viscosity or Forchheimer coefficient, and is supposed to be a scalar, $\mathbf{K}\in\mathbb{R}^{2\times 2}$ is the permeability tensor, and $q$ is a source/sink term. We suppose that $\mathbf{K}$ is a diagonal tensor whose entries $K_{xx}$ and $K_{yy}$ are strictly positive and bounded in $\Omega$. For the sake of convenience, homogeneous Dirichlet boundary conditions are considered, but other types of boundary data can also be handled. We further assume that the porous medium $\Omega$ contains a subset $\Omega_f$ representing a single fracture, which divides the flow domain into two disjoint connected subdomains $\Omega_1$ and $\Omega_2$, i.e.,
$$
\Omega\backslash\overline{\Omega}_f=\Omega_1\cup\Omega_2,\qquad \Omega_1\cap\Omega_2=\emptyset.
$$
In addition, we introduce the notations $\Gamma_{k}=\partial\Omega_k\cap\Gamma$, for $k=1,2,f$, and $\gamma_k=\partial\Omega_k\cap\partial\Omega_f\cap\Omega$, for $k=1,2$. The unit vector normal to $\gamma_k$ pointing outward from $\Omega_k$ is denoted by $\mathbf{n}_k$, for $k=1,2$. A schematic representation of the flow domain including the previous notations is shown in Figure \ref{fig:one:fracture} (left).

Following \cite{fri:rob:saa:08}, we assume that the velocity in the subdomains is small enough to be described by Darcy's law, while that in the fracture needs to be modeled by Forchheimer's law. Under these assumptions, problem \eqref{cont:problem} may be rewritten as the following transmission problem, for $k=1,2,f$, and $j=1,2$
\begin{subequations}\label{transmission:problem}
	\renewcommand{\theequation}{\theparentequation\alph{equation}}
	\begin{align}
	\mathbf{u}_j&=-\mathbf{K}_j\nabla p_j\hspace*{-1.5cm}&&\hbox{in }\Omega_j,\label{transmission:problem:a1}\\
	\left(1+\beta\, |\mathbf{u}_f|\right)\mathbf{u}_f&=-\mathbf{K}_f\nabla p_f\hspace*{-1.5cm}&&\hbox{in }\Omega_f,\label{transmission:problem:a2}\\
	\nabla\cdot\mathbf{u}_k&=q_k\hspace*{-1.5cm}&&\hbox{in }\Omega_k,\label{transmission:problem:b}\\
	p_j&=p_f\hspace*{-1.5cm}&&\hbox{on }\gamma_j,\label{transmission:problem:c}\\
	\mathbf{u}_j\cdot\mathbf{n}_j&=\mathbf{u}_f\cdot\mathbf{n}_j\hspace*{-1.5cm}&&\hbox{on }\gamma_j,\label{transmission:problem:d}\\
	p_k&=0\hspace*{-1.5cm}&&\hbox{on }\Gamma_{k},\label{transmission:problem:e}
	\end{align}
\end{subequations}
where $p_k$, $\mathbf{u}_k$, $\mathbf{K}_k$ and $q_k$ are the restrictions of 
$p$, $\mathbf{u}$, $\mathbf{K}$ and $q$, respectively, to $\Omega_k$, for $k=1,2,f$.
Equations \eqref{transmission:problem:c} and \eqref{transmission:problem:d} provide coupling conditions that guarantee the continuity of the pressure and the normal flux, respectively, across the interfaces between the fracture and the porous matrix.

In the preceding transmission problem, both the bulk and the fracture are defined to be two-dimensional domains. As a consequence, from a numerical viewpoint, we will need extremely fine meshing to resolve the width of the fracture, assumed to be much smaller than its length. This fact will thus increase the computational cost of the algorithm. In order to circumvent this drawback, the fracture is considered to be a one-dimensional interface between the bulk subdomains $\Omega_1$ and $\Omega_2$. The resulting model is known as mixed-dimensional or reduced model. This idea was first proposed in \cite{Martin-Jaffre-Roberts} for a Darcy--Darcy coupling between the fracture and the porous matrix, and has been subsequently used in \cite{kna:rob:14,fri:rob:saa:06,fri:rob:saa:08} in the context of Darcy--Forchheimer couplings. Note that, as an additional advantage in this latter case, the nonlinear Forchheimer problem \eqref{transmission:problem:a2} posed in the fracture is no longer a two-dimensional problem, but a one-dimensional one.


According to \cite{Martin-Jaffre-Roberts}, there exists a non-self-intersecting one-dimensional manifold $\gamma$ such that the fracture can be expressed as
$$
\Omega_f=\left\{\mathbf{x}\in\Omega:\mathbf{x}=\mathbf{s}+\theta\,\mathbf{n}_{\gamma}, \hbox{ for some } \mathbf{s}\in\gamma \hbox{ and } |\theta|<\textstyle\dfrac{d(\mathbf{s})}{2}\right\},
$$
where $d(\mathbf{s})>0$ denotes the width of the fracture at $\mathbf{s}$ in the normal direction, and $\mathbf{n}_{\gamma}$ is the outward unit normal to $\gamma$ with a fixed orientation from $\Omega_1$ to $\Omega_2$. Note that, with this definition, $\mathbf{n}_{\gamma}=\mathbf{n}_1=-\mathbf{n}_2$ (see Figure \ref{fig:one:fracture}). We will suppose that $d(\mathbf{s})$ is much smaller than the other characteristic dimensions of the fracture.

The key point in this procedure is to collapse the fracture $\Omega_f$ into the line $\gamma$, and integrate the equations \eqref{transmission:problem:a2} and \eqref{transmission:problem:b} (the latter for the index $k=f$) along the fracture width. In doing so, we need to split up such equations into their normal and tangential parts. Let us denote the projection operators onto the normal and tangent spaces of $\gamma$ as $\mathbf{P}_{\mathbf{n}}=\mathbf{n}_{\gamma}\mathbf{n}_{\gamma}^T$ and $\mathbf{P}_{\boldsymbol{\tau}}=\mathbf{I}-\mathbf{P}_{\mathbf{n}}$, $\mathbf{I}$ being the identity tensor. For regular vector- and scalar-valued functions $\mathbf{g}$ and $g$, the tangential divergence and gradient operators on the fracture are defined, respectively, as
$$
\nabla^{\boldsymbol{\tau}}\cdot\mathbf{g}=\mathbf{P}_{\boldsymbol{\tau}}:\mathbf{\nabla}\mathbf{g},\qquad\nabla^{\boldsymbol{\tau}} g=\mathbf{P}_{\boldsymbol{\tau}}\nabla g.
$$
Following \cite{for:fum:sco:ruf:14}, we assume that the permeability tensor $\mathbf{K}_f$ decomposes additively as
\begin{equation}\label{K:decomposition}
\mathbf{K}_f=K_{f}^{\mathbf{n}}\mathbf{P}_{\mathbf{n}}+K_{f}^{\boldsymbol{\tau}}\mathbf{P}_{\boldsymbol{\tau}},
\end{equation}
where $K_{f}^{\mathbf{n}}$ and $K_{f}^{\boldsymbol{\tau}}$ are defined to be strictly positive and bounded in $\Omega_f$. Accordingly, $\mathbf{u}_f=\mathbf{u}_{f,\boldsymbol{\tau}}+\mathbf{u}_{f,\mathbf{n}}$, where $\mathbf{u}_{f,\boldsymbol{\tau}}=\mathbf{P}_{\boldsymbol{\tau}}\mathbf{u}_f$ and $\mathbf{u}_{f,\mathbf{n}}=\mathbf{P}_{\mathbf{n}}\mathbf{u}_f$.

In this framework, we introduce the so-called reduced variables, namely: the reduced pressure $p_{\gamma}$, the reduced Darcy velocity $\mathbf{u}_{\gamma}$, and the reduced source/sink term $q_{\gamma}$, formally defined as \cite{Martin-Jaffre-Roberts,delpra:fum:sco:17}
$$
p_{\gamma}(\mathbf{s})=\textstyle\dfrac{1}{d(\mathbf{s})}(p_f,1)_{\ell(\mathbf{s})},\quad\
\mathbf{u}_{\gamma}(\mathbf{s})=(\mathbf{u}_{f,\boldsymbol{\tau}},1)_{\ell(\mathbf{s})},\quad\
q_{\gamma}(\mathbf{s})=(q_f,1)_{\ell(\mathbf{s})},
$$
where $\ell(\mathbf{s})=\left(-\frac{d(\mathbf{s})}{2},\frac{d(\mathbf{s})}{2}\right)$. 
In addition, along the lines of \cite{fri:rob:saa:08}, we assume that the flow in the normal direction within the fracture is described by Darcy's law. This assumption is based on the fact that the ratio between the width and the length of the fracture is small. Thus, equation \eqref{transmission:problem:a2} may be decomposed into its tangential and normal direction as follows
\begin{subequations}\label{Forch:decomp}
	\renewcommand{\theequation}{\theparentequation\alph{equation}}
	\begin{align}
	\left(1+\beta\, |\mathbf{u}_f|\right)\mathbf{u}_{f,\boldsymbol{\tau}}&=-{K}_{f}^{\boldsymbol{\tau}}\,\nabla^{\boldsymbol{\tau}} p_f,\label{Forch:decomp:t}\\
	\mathbf{u}_{f,\mathbf{n}}&=-{K}_{f}^{\mathbf{n}}\,\nabla^{\mathbf{n}} p_f.\label{Forch:decomp:n}
	\end{align}
\end{subequations}

Since $\mathbf{u}_{f,\mathbf{n}}$ is assumed to be much smaller than $\mathbf{u}_{f,\boldsymbol{\tau}}$, we have the approximation $|\mathbf{u}_f|\approx|\mathbf{u}_{f,\boldsymbol{\tau}}|\approx\frac{1}{d}\,|\mathbf{u}_{\gamma}|$.
Then, the integration of \eqref{Forch:decomp:t} along the line segment $\ell(\mathbf{s})$ permits us to derive a Forchheimer's law in the one-dimensional domain $\gamma$. In turn, the integration of \eqref{Forch:decomp:n} in the normal direction to the fracture can be used to give boundary conditions along $\gamma$ for the systems in $\Omega_1$ and $\Omega_2$. Hence, we obtain the following interface problem, for $k=1,2$,
\begin{subequations}
	\renewcommand{\theequation}{\theparentequation\alph{equation}}\label{interface:problem}
	\begin{align}
	\mathbf{u}_k&=-\mathbf{K}_k\nabla p_k&&\hbox{in }\Omega_k,\label{interface:problem:a}\\
	\nabla\cdot\mathbf{u}_k&=q_k&&\hbox{in }\Omega_k,\label{interface:problem:b}\\
	\left(1+\frac{\beta}{d}\,|\mathbf{u}_{\gamma}|\right)\mathbf{u}_{\gamma}&=-dK_{f}^{\boldsymbol{\tau}}\nabla^{\boldsymbol{\tau}} p_{\gamma}&&\hbox{on }\gamma,\label{interface:problem:c}\\
	\nabla^{\boldsymbol{\tau}}\cdot\mathbf{u}_{\gamma}&=q_{\gamma}+(\mathbf{u}_1\cdot\mathbf{n}_1+\mathbf{u}_2\cdot\mathbf{n}_2)&&\hbox{on }\gamma,\label{interface:problem:d}\\
	\alpha_{\gamma}(p_k-p_{\gamma})&=\xi\,\mathbf{u}_k\cdot\mathbf{n}_{k}-(1-\xi)\,\mathbf{u}_{k+1}\cdot\mathbf{n}_{k+1}&&\hbox{on }\gamma,\label{interface:problem:e}\\
	p_k&=0&&\hbox{on }\Gamma_k,\label{interface:problem:f}\\
	p_{\gamma}&=0&&\hbox{on }\partial\gamma,\label{interface:problem:g}
	\end{align}
\end{subequations}
where $\alpha_{\gamma}=2K_{f}^{\mathbf{n}}/d$ and the index $k$ is supposed to vary in $\mathbb{Z}/2\mathbb{Z}$, so that, if $k=2$, then $k+1=1$. According to \cite{Martin-Jaffre-Roberts,ang:boy:hub:09}, $\xi\in(1/2,1]$ is a closure parameter related to the pressure cross profile in the fracture. The ratio $K_{f}^{\mathbf{n}}/d$ and the product $K_{f}^{\boldsymbol{\tau}}d$ are sometimes referred to as effective permeabilities in the normal and tangential directions to the fracture, respectively \cite{fle:fum:sco:16}.

In the preceding system, \eqref{interface:problem:c} represents Forchheimer's law in the tangential direction to the fracture, while \eqref{interface:problem:d} models mass conservation inside the fracture. Remarkably, the additional source term $\mathbf{u}_1\cdot\mathbf{n}_1+\mathbf{u}_2\cdot\mathbf{n}_2$ is introduced on $\gamma$ to take into account the contribution of the subdomain flows to the fracture flow. In turn, \eqref{interface:problem:e} is obtained by averaging the equation \eqref{Forch:decomp:n} in the normal direction to the fracture and using a quadrature rule with weights $\xi$ and $1-\xi$ for integrating $\mathbf{u}_f\cdot\mathbf{n}_k$ across the fracture, for $k=1,2$. Formally, it can be regarded as a Robin boundary condition for the subdomain $\Omega_k$ that involves the pressure in the fracture $p_{\gamma}$ and the normal flux from the neighboring subdomain $\Omega_{k+1}$. It is quite usual to express \eqref{interface:problem:e} in terms of average operators for the pressures and normal fluxes, and jump operators for the pressures across the fracture \cite{delpra:fum:sco:17,ang:sco:12}.

\section{The spatial discretization}\label{sec:discretization}

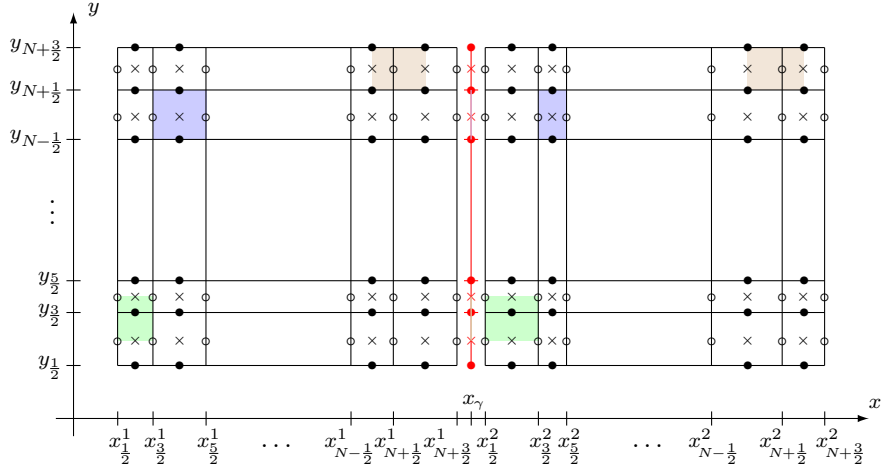
\begin{figure}[t]
		\begin{center}
			\hspace*{-1.1cm}
			\unitlength=0.047cm
			\begin{picture}(200,100)
			\put(-5,5){\vector(1,0){230}}
			\put(0,0){\vector(0,1){120}}
			\put(225,8){{\footnotesize$x$}}
			\put(4,120){{\footnotesize$y$}}
			\put(12.5,3){\line(0,1){4}}
			\put(22.5,3){\line(0,1){4}}
			\put(37.5,3){\line(0,1){4}}
			\put(78.5,3){\line(0,1){4}}
			\put(90.5,3){\line(0,1){4}}
			\put(108.5,3){\line(0,1){4}}
			\put(112.5,3){\line(0,1){4}}
			\put(116.5,3){\line(0,1){4}}
			\put(131.5,3){\line(0,1){4}}
			\put(139.5,3){\line(0,1){4}}
			\put(180.5,3){\line(0,1){4}}
			\put(200.5,3){\line(0,1){4}}
			\put(212.5,3){\line(0,1){4}}
			\put(10,-3){{\footnotesize$x_{\mbox{\hspace*{-0.05cm}\tiny$\frac{1}{2}$}}^1$}}
			\put(20,-3){{\footnotesize$x_{\mbox{\hspace*{-0.05cm}\tiny$\frac{3}{2}$}}^1$}}
			\put(35,-3){{\footnotesize$x_{\mbox{\hspace*{-0.05cm}\tiny$\frac{5}{2}$}}^1$}}
			\put(53,-3){\ldots}
			\put(71,-3){{\footnotesize$x_{\mbox{\hspace*{-0.05cm}\tiny$N\hspace*{-0.08cm}-\hspace*{-0.08cm}\frac{1}{2}$}}^1$}}
			\put(85,-3){{\footnotesize$x_{\mbox{\hspace*{-0.05cm}\tiny$N\hspace*{-0.08cm}+\hspace*{-0.08cm}\frac{1}{2}$}}^1$}}
			\put(99,-3){{\footnotesize$x_{\mbox{\hspace*{-0.05cm}\tiny$N\hspace*{-0.08cm}+\hspace*{-0.08cm}\frac{3}{2}$}}^1$}}
			
			\put(110,9){{\footnotesize$x_{\mbox{\hspace*{-0.02cm}\tiny$\gamma$}}$}}
			
			\put(114,-3){{\footnotesize$x_{\mbox{\hspace*{-0.05cm}\tiny$\frac{1}{2}$}}^2$}}
			\put(129,-3){{\footnotesize$x_{\mbox{\hspace*{-0.05cm}\tiny$\frac{3}{2}$}}^2$}}
			\put(137,-3){{\footnotesize$x_{\mbox{\hspace*{-0.05cm}\tiny$\frac{5}{2}$}}^2$}}
			\put(158,-3){\ldots}   
			\put(174,-3){{\footnotesize$x_{\mbox{\hspace*{-0.05cm}\tiny$N\hspace*{-0.08cm}-\hspace*{-0.08cm}\frac{1}{2}$}}^2$}}
			\put(194,-3){{\footnotesize$x_{\mbox{\hspace*{-0.05cm}\tiny$N\hspace*{-0.08cm}+\hspace*{-0.08cm}\frac{1}{2}$}}^2$}}
			\put(210,-3){{\footnotesize$x_{\mbox{\hspace*{-0.05cm}\tiny$N\hspace*{-0.08cm}+\hspace*{-0.08cm}\frac{3}{2}$}}^2$}}
			
			\put(-2,20){\line(1,0){4}}
			\put(-2,35){\line(1,0){4}}
			\put(-2,44){\line(1,0){4}}
			\put(-2,84){\line(1,0){4}}
			\put(-2,98){\line(1,0){4}}
			\put(-2,110){\line(1,0){4}}
			\put(-10,20){\footnotesize$y_{\mbox{\hspace*{-0.03cm}\tiny$\frac{1}{2}$}}$}
			\put(-10,35){\footnotesize$y_{\mbox{\hspace*{-0.03cm}\tiny$\frac{3}{2}$}}$}
			\put(-10,44){\footnotesize$y_{\mbox{\hspace*{-0.03cm}\tiny$\frac{5}{2}$}}$}
			\put(-7,60){\vdots}
			\put(-18,84){\footnotesize$y_{N\hspace*{-0.03cm}-\mbox{\hspace*{-0.03cm}\tiny$\frac{1}{2}$}}$}
			\put(-18,98){\footnotesize$y_{N\hspace*{-0.03cm}+\mbox{\hspace*{-0.03cm}\tiny$\frac{1}{2}$}}$}
			\put(-18,110){\footnotesize$y_{N\hspace*{-0.03cm}+\mbox{\hspace*{-0.03cm}\tiny$\frac{3}{2}$}}$}
			
			{\color{blue!20}\polygon*(22.5,84)(37.5,84)(37.5,98)(22.5,98)}
			{\color{blue!20}\polygon*(129,84)(137,84)(137,98)(129,98)}
			
			{\color{green!20}\polygon*(7.5,27)(17.5,27)(17.5,39.5)(7.5,39.5)}
			{\color{green!20}\polygon*(109,27)(124,27)(124,39.5)(109,39.5)}

			{\color{brown!20}\polygon*(74.5,98)(89.5,98)(89.5,110)(74.5,110)}
			{\color{brown!20}\polygon*(178,98)(194,98)(194,110)(178,110)}

			\put(0,20){\line(1,0){96}}
			\put(104,20){\line(1,0){96}}\put(200,20){\line(0,1){90}}\put(0,110){\line(1,0){96}}
			\put(104,110){\line(1,0){96}}\put(0,20){\line(0,1){90}}
			\put(100,20){\color{red}\line(0,1){90}}
			\put(96,20){\line(0,1){90}}
			\put(104,20){\line(0,1){90}}
			\put(10,20){\line(0,1){90}}
			\put(25,20){\line(0,1){90}}
			\put(66,20){\line(0,1){90}}
			\put(78,20){\line(0,1){90}}
			\put(0,35){\line(1,0){96}}
			\put(0,44){\line(1,0){96}}
			\put(0,84){\line(1,0){96}}
			\put(0,98){\line(1,0){96}}
			\put(98,35){\color{red}\line(1,0){4}}
			\put(98,44){\color{red}\line(1,0){4}}
			\put(98,84){\color{red}\line(1,0){4}}
			\put(98,98){\color{red}\line(1,0){4}}
			\put(104,35){\line(1,0){96}}
			\put(104,44){\line(1,0){96}}
			\put(104,84){\line(1,0){96}}
			\put(104,98){\line(1,0){96}}
			\put(119,20){\line(0,1){90}}
			\put(127,20){\line(0,1){90}}
			\put(188,20){\line(0,1){90}}
			\put(168,20){\line(0,1){90}}
			\put(3,26){{\tiny$\times$}}
			\put(15.5,26){{\tiny$\times$}}
			\put(70,26){{\tiny$\times$}}
			\put(85,26){{\tiny$\times$}}
			\put(3,38.5){{\tiny$\times$}}
			\put(15.5,38.5){{\tiny$\times$}}
			\put(70,38.5){{\tiny$\times$}}
			\put(85,38.5){{\tiny$\times$}}
			\put(3,89.5){{\tiny$\times$}}
			\put(15.5,89.5){{\tiny$\times$}}
			\put(70,89.5){{\tiny$\times$}}
			\put(85,89.5){{\tiny$\times$}}
			\put(3,103){{\tiny$\times$}}
			\put(15.5,103){{\tiny$\times$}}
			\put(70,103){{\tiny$\times$}}
			\put(85,103){{\tiny$\times$}}
			\put(98,26){{\color{red}\tiny$\times$}}
			\put(98,38.5){{\color{red}\tiny$\times$}}
			\put(98,89.5){{\color{red}\tiny$\times$}}
			\put(98,103){{\color{red}\tiny$\times$}}
			\put(109.5,26){{\tiny$\times$}}
			\put(121,26){{\tiny$\times$}}
			\put(176.2,26){{\tiny$\times$}}
			\put(192,26){{\tiny$\times$}}
			\put(109.5,38.5){{\tiny$\times$}}
			\put(121,38.5){{\tiny$\times$}}
			\put(176.2,38.5){{\tiny$\times$}}
			\put(192,38.5){{\tiny$\times$}}
			\put(109.5,89.5){{\tiny$\times$}}
			\put(121,89.5){{\tiny$\times$}}
			\put(176.2,89.5){{\tiny$\times$}}
			\put(192,89.5){{\tiny$\times$}}
			\put(109.5,103){{\tiny$\times$}}
			\put(121,103){{\tiny$\times$}}
			\put(176.2,103){{\tiny$\times$}}
			\put(192,103){{\tiny$\times$}}
			\put(3.5,18.7){{\scriptsize$\bullet$}}
			\put(16,18.7){{\scriptsize$\bullet$}}
			\put(70.5,18.7){{\scriptsize$\bullet$}}
			\put(85.5,18.7){{\scriptsize$\bullet$}}
			\put(3.5,33.7){{\scriptsize$\bullet$}}
			\put(16,33.7){{\scriptsize$\bullet$}}
			\put(70.5,33.7){{\scriptsize$\bullet$}}
			\put(85.5,33.7){{\scriptsize$\bullet$}}
			\put(3.5,42.7){{\scriptsize$\bullet$}}
			\put(16,42.7){{\scriptsize$\bullet$}}
			\put(70.5,42.7){{\scriptsize$\bullet$}}
			\put(85.5,42.7){{\scriptsize$\bullet$}}
			\put(3.5,82.7){{\scriptsize$\bullet$}}
			\put(16,82.7){{\scriptsize$\bullet$}}
			\put(70.5,82.7){{\scriptsize$\bullet$}}
			\put(85.5,82.7){{\scriptsize$\bullet$}}
			\put(3.5,96.7){{\scriptsize$\bullet$}}
			\put(16,96.7){{\scriptsize$\bullet$}}
			\put(70.5,96.7){{\scriptsize$\bullet$}}
			\put(85.5,96.7){{\scriptsize$\bullet$}}
			\put(3.5,108.7){{\scriptsize$\bullet$}}
			\put(16,108.7){{\scriptsize$\bullet$}}
			\put(70.5,108.7){{\scriptsize$\bullet$}}
			\put(85.5,108.7){{\scriptsize$\bullet$}}
			\put(98.5,18.7){{\color{red}\scriptsize$\bullet$}}
			\put(98.5,33.7){{\color{red}\scriptsize$\bullet$}}
			\put(98.5,42.7){{\color{red}\scriptsize$\bullet$}}
			\put(98.5,82.7){{\color{red}\scriptsize$\bullet$}}
			\put(98.5,96.7){{\color{red}\scriptsize$\bullet$}}
			\put(98.5,108.7){{\color{red}\scriptsize$\bullet$}}
			\put(110,18.7){{\scriptsize$\bullet$}}
			\put(121.5,18.7){{\scriptsize$\bullet$}}
			\put(176.8,18.7){{\scriptsize$\bullet$}}
			\put(192.7,18.7){{\scriptsize$\bullet$}}
			\put(110,33.7){{\scriptsize$\bullet$}}
			\put(121.5,33.7){{\scriptsize$\bullet$}}
			\put(176.8,33.7){{\scriptsize$\bullet$}}
			\put(192.7,33.7){{\scriptsize$\bullet$}}
			\put(110,42.7){{\scriptsize$\bullet$}}
			\put(121.5,42.7){{\scriptsize$\bullet$}}
			\put(176.8,42.7){{\scriptsize$\bullet$}}
			\put(192.7,42.7){{\scriptsize$\bullet$}}
			\put(110,82.7){{\scriptsize$\bullet$}}
			\put(121.5,82.7){{\scriptsize$\bullet$}}
			\put(176.8,82.7){{\scriptsize$\bullet$}}
			\put(192.7,82.7){{\scriptsize$\bullet$}}
			\put(110,96.7){{\scriptsize$\bullet$}}
			\put(121.5,96.7){{\scriptsize$\bullet$}}
			\put(176.8,96.7){{\scriptsize$\bullet$}}
			\put(192.7,96.7){{\scriptsize$\bullet$}}
			\put(110,108.7){{\scriptsize$\bullet$}}
			\put(121.5,108.7){{\scriptsize$\bullet$}}
			\put(176.8,108.7){{\scriptsize$\bullet$}}
			\put(192.7,108.7){{\scriptsize$\bullet$}}
			\put(-1.5,25.5){{\scriptsize$\circ$}}
			\put(8.5,25.5){{\scriptsize$\circ$}}
			\put(23.5,25.5){{\scriptsize$\circ$}}
			\put(64.5,25.5){{\scriptsize$\circ$}}
			\put(76.5,25.5){{\scriptsize$\circ$}}
			\put(94.5,25.5){{\scriptsize$\circ$}}
			\put(-1.5,38){{\scriptsize$\circ$}}
			\put(8.5,38){{\scriptsize$\circ$}}
			\put(23.5,38){{\scriptsize$\circ$}}
			\put(64.5,38){{\scriptsize$\circ$}}
			\put(76.5,38){{\scriptsize$\circ$}}
			\put(94.5,38){{\scriptsize$\circ$}}
			\put(-1.5,89){{\scriptsize$\circ$}}
			\put(8.5,89){{\scriptsize$\circ$}}
			\put(23.5,89){{\scriptsize$\circ$}}
			\put(64.5,89){{\scriptsize$\circ$}}
			\put(76.5,89){{\scriptsize$\circ$}}
			\put(94.5,89){{\scriptsize$\circ$}}
			\put(-1.5,102.5){{\scriptsize$\circ$}}
			\put(8.5,102.5){{\scriptsize$\circ$}}
			\put(23.5,102.5){{\scriptsize$\circ$}}
			\put(64.5,102.5){{\scriptsize$\circ$}}
			\put(76.5,102.5){{\scriptsize$\circ$}}
			\put(94.5,102.5){{\scriptsize$\circ$}}
			\put(102.5,25.5){{\scriptsize$\circ$}}
			\put(117.5,25.5){{\scriptsize$\circ$}}
			\put(125.5,25.5){{\scriptsize$\circ$}}
			\put(166.5,25.5){{\scriptsize$\circ$}}
			\put(186.5,25.5){{\scriptsize$\circ$}}
			\put(198.5,25.5){{\scriptsize$\circ$}}
			\put(102.5,38){{\scriptsize$\circ$}}
			\put(117.5,38){{\scriptsize$\circ$}}
			\put(125.5,38){{\scriptsize$\circ$}}
			\put(166.5,38){{\scriptsize$\circ$}}
			\put(186.5,38){{\scriptsize$\circ$}}
			\put(198.5,38){{\scriptsize$\circ$}}
			\put(102.5,89){{\scriptsize$\circ$}}
			\put(117.5,89){{\scriptsize$\circ$}}
			\put(125.5,89){{\scriptsize$\circ$}}
			\put(166.5,89){{\scriptsize$\circ$}}
			\put(186.5,89){{\scriptsize$\circ$}}
			\put(198.5,89){{\scriptsize$\circ$}}
			\put(102.5,102.5){{\scriptsize$\circ$}}
			\put(117.5,102.5){{\scriptsize$\circ$}}
			\put(125.5,102.5){{\scriptsize$\circ$}}
			\put(166.5,102.5){{\scriptsize$\circ$}}
			\put(186.5,102.5){{\scriptsize$\circ$}}
			\put(198.5,102.5){{\scriptsize$\circ$}}
			
			\put(100,84){\color{blue!20}\line(0,1){14}}
			
			\put(100,27){\color{green!20}\line(0,1){12.5}}

			\end{picture}
		\end{center}\vspace*{0.1cm}
		\caption{Staggered grid location of unknowns and corresponding control volumes.}
		\label{fig:staggered:grid}
	\end{figure}

	Let us assume that the subdomains $\Omega_k$ admit rectangular partitions $\mathcal{T}_{h}^{k}$, for $k=1,2$, that match at the interface $\gamma$.
	Such meshes $\mathcal{T}_{h}^{k}$ induce a unique partition on $\gamma$ denoted by $\mathcal{T}_h^{\gamma}$. In the case of considering a vertical fracture 
	as that shown in Figure \ref{fig:one:fracture}, such partitions may be defined as  $\mathcal{T}_{h}^{k}=\cup_{i,j=1}^{N+1}E_{i,j}^k$ and $\mathcal{T}_h^{\gamma}=\cup_{j=1}^{N+1}E_{j}^{\gamma}$, where
		\begin{align*}E_{i,j}^k&=(x_{i-1/2}^k,x_{i+1/2}^k)\times(y_{j-1/2},y_{j+1/2}), \\
		E_j^{\gamma}&=\{x_{\gamma}\}\times(y_{j-1/2},y_{j+1/2}),
			\end{align*}
	$x_{\gamma}$ being equal to $x_{N+3/2}^1$ and $x_{1/2}^2$. In the framework of finite volume methods, these sets are known as control volumes. Figure \ref{fig:staggered:grid} shows the control volumes $E_{2,N}^1$, $E_{2,N}^2$ and $E_N^\gamma$ highlighted in blue. Note that both $E_{2,N}^1$ and $E_{2,N}^2$ are two-dimensional control volumes, while $E_N^\gamma$ is one-dimensional.
	
	In this setting, we associate the pressure unknowns $p_{i,j}^k$ and $p_j^{\gamma}$ to the element centers $(x_i^k,y_j)$ and $(x_\gamma,y_j)$, respectively, as indicated with cross signs in Figure \ref{fig:staggered:grid}. In particular,
	$$
	p_{i,j}^k\approx \frac{1}{|E_{i,j}^k|}\iint_{E_{i,j}^k} p(x,y)\,dx\,dy,\qquad p_{j}^{\gamma}\approx \frac{1}{|E_{j}^\gamma|}\int_{E_{j}^\gamma} p_{\gamma}(s)\,ds.
	$$
	In order to introduce the velocity unknowns, we first define some additional control volumes associated to the midpoints of the edges of the meshes $\mathcal{T}_h^k$, for $k=1,2$, and $\mathcal{T}_h^{\gamma}$. In particular, let us define the following control volumes associated to the vertical edges of $\mathcal{T}_h^k$,
	\begin{align*}
	E_{i+1/2,j}^k&=(x_{i}^k,x_{i+1}^k)\times(y_{j-1/2},y_{j+1/2}),\qquad i=1,\ldots,N,\ j=1,\ldots,N+1,\\
	E_{1/2,j}^k&=(x_{1/2}^k,x_{1}^k)\times(y_{j-1/2},y_{j+1/2}),\qquad j=1,\ldots,N+1,\\
	E_{N+3/2,j}^k&=(x_{N+1}^k,x_{N+3/2}^k)\times(y_{j-1/2},y_{j+1/2}),\qquad j=1,\ldots,N+1. 
	\end{align*}
	The control volumes associated to the horizontal edges of $\mathcal{T}_h^k$ are denoted by $E_{i,j+1/2}^k$, for $i=1,\ldots,N+1$ and $j=0,1,\ldots,N+1$, and may be defined in a similar way.
	Finally, we define the following one-dimensional control volumes associated to the  mesh points $(x_\gamma,y_{j+1/2})$ of the partition $\mathcal{T}_h^\gamma$,
	\begin{align*}
	E_{j+1/2}^\gamma&=\{x_{\gamma}\}\times(y_{j},y_{j+1}),\qquad j=1,\ldots,N,\\
	E_{1/2}^\gamma&=\{x_{\gamma}\}\times(y_{1/2},y_{1}),\\
	E_{N+3/2}^\gamma&=\{x_{\gamma}\}\times(y_{N+1},y_{N+3/2}).
	\end{align*}
	Figure \ref{fig:staggered:grid} shows the control volumes $E_{N+1/2,N+1}^1$ and $E_{N+1/2,N+1}^2$ highlighted in brown, and the  control volumes $E_{1,3/2}^1$, $E_{1,3/2}^2$ and $E_{3/2}^\gamma$ highlighted in green.
	
	In this context, the velocity unknowns can be classified into three groups. The first group comprises the normal flux components associated to the vertical edges of the two-dimen\-sional grids $\mathcal{T}_{h}^{k}$, which are denoted by $u_{i+1/2,j}^k$, for $i=0,1,\ldots,N+1$, $j=1,\ldots,N+1$, $k=1,2$, and are represented by black empty circles in Figure \ref{fig:staggered:grid}. The second set contains the normal flux components associated to the horizontal edges of $\mathcal{T}_{h}^{k}$, which are denoted by $v_{i,j+1/2}^k$, for $i=1,\ldots,N+1$, $j=0,1,\ldots,N+1$, $k=1,2$, and are depicted by black filled dots in the same plot. Finally, the third group comprises the normal flux components associated to the edges of the one-dimensional grid $\mathcal{T}_{h}^{\gamma}$, which are denoted by $u_{j+1/2}^\gamma$, for $j=0,1,\ldots,N+1$, and are marked by red filled dots.	
	In particular, 
	\begin{align*}
	u_{i+1/2,j}^k&\approx \frac{1}{|E_{i+1/2,j}^k|}\iint_{E_{i+1/2,j}^k}u^k(x,y)\,dx\,dy,\\[1ex]
	v_{i,j+1/2}^k&\approx \frac{1}{|E_{i,j+1/2}^k|}\iint_{E_{i,j+1/2}^k}v^k(x,y)\,dx\,dy,\\[1ex]
	u_{j+1/2}^{\gamma}&\approx	\frac{1}{|E_{j+1/2}^{\gamma}|}\int_{E_{j+1/2}^{\gamma}}u^{\gamma}(s)\,ds.
	\end{align*}
Let us introduce the notation $\mathbf{u}^k=(u^k,v^k)^T$ for the two components of the velocity on $\Omega_k$, for $k=1,2$. Taking into account that $\mathbf{K}_k$ are diagonal tensors, with diagonal coefficients $K_{xx}^k$ and $K_{yy}^k$, the equation (\ref{interface:problem:a}) can be decomposed as
\begin{subequations}\label{decomp:darcy:eq}
	\renewcommand{\theequation}{\theparentequation\alph{equation}}
	\begin{align}
		u^k+K_{xx}^k\,\frac{\partial p_k}{\partial x}&=0,\label{decomp:darcy:eq:x}\\
		v^k+K_{yy}^k\,\frac{\partial p_k}{\partial y}&=0.\label{decomp:darcy:eq:y}
	\end{align}
\end{subequations}
The integration of equation \eqref{decomp:darcy:eq:x} over the control volumes $E_{i-1/2,j}^k$, together with the application of the midpoint quadrature rule in the $x$-direction and a suitable approximation of the flux at the midpoints of the cell edges, gives rise to the discrete equations for the horizontal velocities. In this case, we consider that the flux over each edge is approximated by using the pressure unknowns in the two cells sharing that edge. This scheme, known as the two-point flux approximation method \cite{eym:gal:gui:her:mas:14}, is widely used in reservoir simulations. In turn, the discrete equations for the vertical velocities are obtained by integrating equation \eqref{decomp:darcy:eq:y} over the control volumes $E_{i,j-1/2}$ and following a similar procedure.

In particular, the interior velocity unknowns will satisfy the following equations
\begin{align*}
u_{i-1/2,j}^k+2\left(\frac{\Delta x_{i-1}^k}{(K_{xx}^k)_{i-1,j}}+\frac{\Delta x_{i}^k}{(K_{xx}^k)_{i,j}}\right)^{-1}\,(p_{i,j}^k-p_{i-1,j}^k)&=0,
\\
v_{i,j-1/2}^k+2\left(\frac{\Delta y_{j-1}}{(K_{yy}^k)_{i,j-1}}+\frac{\Delta y_j}{(K_{yy}^k)_{i,j}}\right)^{-1}\,(p_{i,j}^k-p_{i,j-1}^k)&=0,
\end{align*}
where $\Delta x_i^k=x_{i+1/2}^k-x_{i-1/2}^k$ and $\Delta y_j=y_{j+1/2}-y_{j-1/2}$.
Finally, integrating equation \eqref{interface:problem:b} over the control volumes $E_{i,j}^k$ and applying the divergence theorem, we get
$$
\frac{u_{i+1/2,j}^k-u_{i-1/2,j}^k}{\Delta x_i^k}+\frac{v_{i,j+1/2}^k-v_{i,j-1/2}^k}{\Delta y_j}=q_{i,j}^k, 
$$
where $$q_{i,j}^{k}=\frac{1}{\Delta x_i^k\,\Delta y_j}\iint_{E_{i,j}^k}q_k\,dx\,dy.$$ 
Similarly, by integrating equations \eqref{interface:problem:c} and \eqref{interface:problem:d} over the control volumes $E_{j-1/2}^\gamma$ and $E_{j}^\gamma$, respectively, we get the following equations for the interior unknowns of $\gamma$
\begin{align*}
\left(1+\displaystyle\frac{\beta}{d}|u_{j-1/2}^{\gamma}|\right)u_{j-1/2}^{\gamma}+2\,d\left(\displaystyle\frac{\Delta y_{j-1}}{(K_{f}^{\boldsymbol{\tau}})_{j-1}}+\frac{\Delta y_j}{(K_{f}^{\boldsymbol{\tau}})_{j}}\right)^{-1}\,(p_{j}^{\gamma}-p_{j-1}^{\gamma})&=0,\\[2ex]
\displaystyle\frac{u_{j+1/2}^{\gamma}-u_{j-1/2}^{\gamma}}{\Delta y_j}-(u_{N+3/2,j}^1-u_{1/2,j}^2)&=q_{j}^{\gamma}, 
\end{align*}
where $q_j^\gamma=\frac{1}{|E_j^\gamma|}\int_{E_j^\gamma}q_\gamma(s)\,ds.$

Considering homogeneous Dirichlet boundary conditions, the equations for the normal fluxes at the horizontal boundaries are given by
\begin{align*}
v_{i,1/2}^k+2\,\displaystyle\frac{(K_{yy}^k)_{i,1}}{\Delta y_1}\,p_{i,1}^k&=0, 
\\[2ex]
v_{i,N+3/2}^k-2\,\displaystyle\frac{(K_{yy}^k)_{i,N+1}}{\Delta y_{N+1}}\,p_{i,N+1}^k&=0,
\\[2ex]
\left(1+\displaystyle\frac{\beta}{d}|u_{1/2}^{\gamma}|\right)u_{1/2}^{\gamma}+2\,d\,\displaystyle\frac{(K_{f}^{\boldsymbol{\tau}})_1\,}{\Delta y_1}\,p_1^{\gamma}&=0,
\\[2ex]
\left(1+\displaystyle\frac{\beta}{d}|u_{N+3/2}^{\gamma}|\right)u_{N+3/2}^{\gamma}-2\,d\,\displaystyle\frac{(K_{f}^{\boldsymbol{\tau}})_{N+1}\,}{\Delta y_{N+1}}\,p_{N+1}^{\gamma}&=0.
\end{align*}
In turn, the equations for the normal fluxes at the vertical boundaries take the form
\begin{align*}
u_{1/2,j}^1+2\,\displaystyle\frac{(K_{xx}^1)_{1,j}}{\Delta x_1^1}\,p_{1,j}^1&=0,\\[2ex]
u_{N+3/2,j}^2-2\,\displaystyle\frac{(K_{xx}^2)_{N+1,j}}{\Delta x_{N+1}^2}\,p_{N+1,j}^2&=0.
\end{align*}
Finally, considering the coupling condition \eqref{interface:problem:e}, the equations for the normal fluxes of the porous matrix at the interface $\gamma$ are
\begin{align*}
\left(1+\frac{\xi\, a_j^1}{\alpha_{\gamma}}\right)u_{N+3/2,j}^1+a_j^1\,(p_j^{\gamma}-p_{N+1,j}^1)+\frac{(1-\xi)\, a_j^1}{\alpha_{\gamma}}\, u_{1/2,j}^2&=0,\\[2ex] 
\left(1+\frac{\xi\,a_j^2}{\alpha_{\gamma}}\right)u_{1/2,j}^2+a_j^2\,(p_{1,j}^2-p_j^{\gamma})+\frac{(1-\xi)\,a_j^2}{\alpha_{\gamma}}\, u_{N+3/2,j}^1&=0, 
\end{align*}
for $j=1,\ldots,N+1$, where $a^1_j=2(K_{xx}^1)_{N+1,j}/\Delta x_{N+1}^1$ and $a_j^2=2 (K_{xx}^2)_{1,j}/\Delta x_1^2$.

Suitable scaling of the previous equations results in a nonlinear saddle point problem of the form
	$$\begin{pmatrix}
	A_1 & C^T & 0 & B_1^T & 0 & F_1^T\\[0.5ex]
	C & A_2 & 0 & 0 & B_2^T & F_2^T\\[0.5ex]
	0 & 0 & A_{\gamma}(U_{\gamma}) & 0 & 0 & B_{\gamma}^T\\[0.5ex]
	B_1 & 0 & 0 & 0 & 0 & 0\\[0.5ex]
	0 & B_2 & 0 & 0 & 0 & 0\\[0.5ex]
	F_1 & F_2 & B_{\gamma} & 0 & 0 & 0
	\end{pmatrix} 
	\begin{pmatrix}
	U_1\\[0.5ex]U_2\\[0.5ex]U_{\gamma}\\[0.5ex]P_1\\[0.5ex]P_2\\[0.5ex]P_{\gamma}
	\end{pmatrix}=
	\begin{pmatrix}
	0\\[0.5ex]0\\[0.5ex]0\\[0.5ex]Q_1\\[0.5ex]Q_2\\[0.5ex]Q_{\gamma}
	\end{pmatrix},
	$$
	where $U^1$, $U^2$ and $U^\gamma$ comprise the velocity unknowns on $\Omega_1$, $\Omega_2$ and $\gamma$, respectively. Similarly, $P^1$, $P^2$ and $P^\gamma$ contain the pressure unknowns on $\Omega_1$, $\Omega_2$ and $\gamma$, respectively. The matrices $A_1$, $A_2$ and $A_{\gamma}(U_{\gamma})$ are diagonal. It is important to notice that the nonlinearity of the problem is associated only with the one-dimensional Forchheimer equation posed on the fracture.

\section{The monolithic multigrid method}\label{sec:multigrid}

Typically, there are two approaches for solving nonlinear problems by using multigrid techniques. One is to apply some linearization method, such as Newton's iteration or Picard method, and then to use multigrid for solving the linear problem corresponding to each iteration step. The second approach, known as full approximation scheme (FAS) \cite{Bra77} consists of applying multigrid directly to the nonlinear problem. This is the approach followed in this work, and it is briefly described here.  It is only needed to describe the two-level version since, similar to the linear case, the nonlinear FAS multigrid method can be defined recursively on the basis the two-level method. In this way, if $A_h(u_h) = f_h$ denotes a nonlinear system of the equations, a two-level FAS scheme reads as follows\\

\noindent {\sf Full Approximation Scheme (FAS)}:
\begin{itemize}
	\item Pre-smoothing: Compute $\bar{u}_h^m$ by applying $\nu_1$ smoothing steps: $\bar{u}_h^m = S_h^{\nu_1} u_h^m$.
	\item Restrict the residual and the current approximation to the coarse grid: $r_H = I_{h,H} (f_h - A_h(\bar{u}_h^m))$, and $u_H^m = \widetilde{I}_{h,H} \bar{u}_h^m$.
	\item Solve the coarse-grid problem $A_H(v_H^m) = A_H(u_H^m) + r_H$.
	\item Interpolate the error approximation to the fine grid and correct the current fine grid approximation: $\hat{u}_h^m = \bar{u}_h^m + I_{H,h} (v_H^m - u_H^m)$.
	\item Post-smoothing: Compute $u_h^{m+1}$ by applying $\nu_2$ smoothing steps : $u_h^{m+1} = S_h^{\nu_2} \hat{u}_h^m$.
\end{itemize}
Here $S_h$ denotes a nonlinear relaxation procedure, $I_{h,H}, \widetilde{I}_{h,H}$ are, possibly different, transfer operators from the fine to the coarse grid, and $I_{H,h}$ is a transfer operator from the coarse to the fine grid. Notice that if $A_h$ is a linear operator, then the FAS scheme is identical to the standard linear multigrid method. 

In this work, we propose a monolithic multigrid method for solving the Darcy--Forchheimer flow in a fractured porous media. This means that we do not iterate between the subproblems in the matrix and within the fracture network, and we treat the whole problem at once. It is well known that the performance of a multigrid method strongly depends on its components. Notice that we are considering a nonlinear saddle point problem with a mixed-dimensional character, and this will have a great influence on the choice of the multigrid elements, which will have to appropriately deal with these characteristics. Next, we describe the components used to define the nonlinear multigrid scheme.

\subsection{Inter-grid transfer operators}
In this section, we introduce the restriction and interpolation operators involved in the multigrid method for solving the mixed-dimensional problem. We consider different transfer operators for the unknowns belonging to the matrix
and for those located at the fractures. In particular, we choose two-dimensional and one-dimensional transfer operators, respectively. This means that we implement mixed-dimensional transfer operators in our multigrid algorithm in order to handle the problem at once. Regarding the unknowns of the porous medium, we take into account the staggered arrangement of their location. Thus, the inter-grid transfer operators that act in the 
porous media unknowns are defined as follows: a six-point restriction is considered at velocity grid points, 
and a four-point restriction is applied at pressure grid points. In stencil notation, these restriction operators are given by
$$
I_{h,H}^{u}=\frac{1}{8}\begin{pmatrix}
1 & 2 & 1\\  &* & \\ 1 & 2 & 1
\end{pmatrix}_h, \quad
I_{h,H}^{v}=\frac{1}{8}\begin{pmatrix}
1 & & 1\\ 2 &* & 2\\ 1 & & 1
\end{pmatrix}_h, \quad
I_{h,H}^{p}=\frac{1}{4}\begin{pmatrix}
1 &  & 1\\  &*&  \\ 1 & & 1
\end{pmatrix}_h,
$$
respectively. We have used the same restriction operators for the current approximations, that is, $\widetilde{I}_{h,H}^u = I_{h,H}^{u}$,  $\widetilde{I}_{h,H}^v = I_{h,H}^{v}$ and  $\widetilde{I}_{h,H}^p = I_{h,H}^{p}$. As the prolongation operators $I_{H,h}^{u/v/p}$, we choose the adjoints of the restrictions.

Regarding the inter-grid transfer operators for the unknowns at the fractures, we again take into account their one-dimensional staggered arrangement, yielding the following restriction transfer operators
$$
\widetilde{I}_{h,H}^{u^{\gamma}} = I_{h,H}^{u^{\gamma}}=\frac{1}{4}\begin{pmatrix}
1 & 2 & 1
\end{pmatrix}_h, \qquad
\widetilde{I}_{h,H}^{p^{\gamma}} = I_{h,H}^{p^{\gamma}}=\frac{1}{2}\begin{pmatrix}
1 & * & 1
\end{pmatrix}_h.
$$
On the other hand, the prolongation operators are their adjoints.
\subsection{Smoother}\label{sec:smoother}
The proposed smoother is based on the well-known Vanka relaxation procedure, which was proposed by Vanka in \cite{vanka} for solving the staggered finite difference discretization of the Navier--Stokes equations. This smoother was based on simultaneously updating all unknowns appearing in the discrete divergence operator in the pressure equation. 
Thus, the relaxation designed here for the proposed coupled problem is based on combining two- and one-dimensional Vanka-type smoothers for the unknowns in the matrix and within the fractures, respectively. Moreover, in each smoothing step two relaxations of the one-dimensional Vanka smoother with a relaxation parameter $\omega = 0.7$ are performed within the fractures, whereas only one smoothing iteration of the two-dimensional Vanka relaxation is carried out within the matrix. Notice that the computational cost of the smoothing step in the fracture is negligible in comparison with that of the porous matrix, since it is a one-dimensional calculation.

\begin{figure}[t]
	\begin{center}
		\unitlength=0.047cm
		\begin{picture}(260,94)
		\put(30,25){\line(1,0){60}}\put(85,20){\line(0,1){60}}
		\put(30,75){\line(1,0){60}}\put(35,20){\line(0,1){60}}
		\dashline[20]{2}(50,20)(50,40)
		\dashline[20]{2}(30,40)(50,40)
		\dashline[20]{2}(30,40)(30,60)
		\dashline[20]{2}(30,60)(50,60)
		\dashline[20]{2}(50,60)(50,80)
		\dashline[20]{2}(50,80)(70,80)
		\dashline[20]{2}(70,80)(70,60)
		\dashline[20]{2}(70,60)(90,60)
		\dashline[20]{2}(90,60)(90,40)
		\dashline[20]{2}(70,40)(90,40)
		\dashline[20]{2}(70,20)(70,40)
		\dashline[20]{2}(50,20)(70,20)
		
		\put(56.5,47){\Large$\times$}
		\put(32.5,47){\Large$\circ$}
		\put(82.5,47){\Large$\circ$}
		\put(58,22.5){\Large$\bullet$}
		\put(58,72.5){\Large$\bullet$}
								
		\put(56.5,42){$p_{i,j}$}
		\put(3.5,47){{\color{blue}$u_{i-1/2,j}$}}
		\put(92,47){{\color{blue}$u_{i+1/2,j}$}}
		\put(53,14){{\color{blue}$v_{i,j-1/2}$}}
		\put(53,85){{\color{blue}$v_{i,j+1/2}$}}

		\put(140,25){\line(1,0){100}}\put(235,20){\line(0,1){60}}
		\put(140,75){\line(1,0){100}}\put(145,20){\line(0,1){60}}
		\put(190,20){\line(0,1){60}}
		\dashline[20]{2}(179,20)(179,80)
		\dashline[20]{2}(179,80)(201,80)
		\dashline[20]{2}(179,20)(201,20)
		\dashline[20]{2}(201,20)(201,80)
		
		\put(186,47){\Large$\times$}
		\put(187.5,22.5){\Large$\bullet$}
		\put(187.5,72.3){\Large$\bullet$}
		
		\put(193.5,47){$p_j^\gamma$}
		\put(190,12){{\color{blue}$u_{j-1/2}^\gamma$}}
		\put(190,85){{\color{blue}$u_{j+1/2}^\gamma$}}
		
		\put(148,30){{\color{red}$E_{N+1,j}^1$}}
		\put(219,30){{\color{red}$E_{1,j}^2$}}
	
		\end{picture}
	\end{center}\vspace*{-1cm}
	\caption{Unknowns updated together by the (left) two-dimensional and (right) one-dimensional Vanka-type smoothers, applied in the porous matrix and within the fractures, respectively.}
	\label{vanka_blocks}
\end{figure}
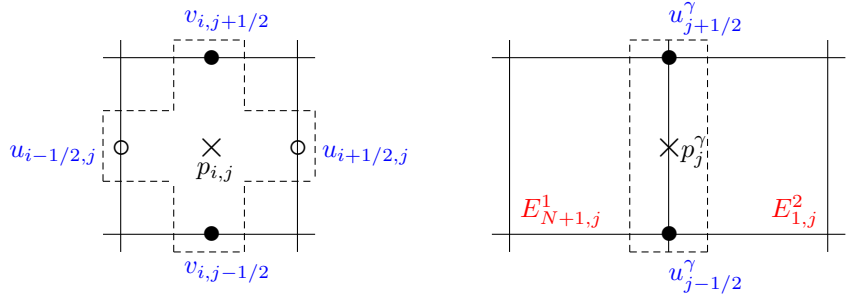


More concretely, in the case of the porous matrix, the Vanka smoothing approach implies that four unknowns corresponding to velocities and one pressure unknown are simultaneously updated, see Figure \ref{vanka_blocks} (left). This means that within the two-dimensional smoothing step, we iterate over all cells and for each cell a $5 \times 5$ system is solved. In particular, the system to solve on each cell, in terms of increments, is written as follows,
$$
\begin{pmatrix}
1 & 0 & 0 & 0 & -\frac{K_{xx}}{h}  \\[0.5ex]
0 & 1 & 0 & 0 & \frac{K_{xx}}{h} \\[0.5ex]
0 & 0 & 1 & 0 & -\frac{K_{yy}}{h}  \\[0.5ex]
0 & 0 & 0 & 1 & \frac{K_{yy}}{h} \\[0.5ex]
\frac{1}{h} & -\frac{1}{h} & \frac{1}{h} & -\frac{1}{h} & 0
\end{pmatrix}
\begin{pmatrix}
\delta u_{i+1/2,j} \\[0.5ex]
\delta u_{i-1/2,j} \\[0.5ex]
\delta v_{i,j+1/2} \\[0.5ex]
\delta v_{i,j-1/2} \\[0.5ex]
\delta p_{i,j} 
\end{pmatrix} = 
\begin{pmatrix}
r^u_{i+1/2,j} \\[0.5ex]
r^u_{i-1/2,j} \\[0.5ex]
r^v_{i,j+1/2} \\[0.5ex]
r^v_{i,j-1/2} \\[0.5ex]
r^p_{i,j} 
\end{pmatrix}.
$$
For the sake of simplicity, in the previous system we omitted the superscript $k$ indicating the corresponding subdomain, we considered a uniform grid in both directions with mesh-size $h$ and we assumed that the permeability tensors are homogeneous on the corresponding subdomain. 

Regarding the one-dimensional Vanka smoother for the fractures, three unknowns are simultaneously updated for each cell: the unknown corresponding to the pressure and both velocities included in the cell, see Figure \ref{vanka_blocks} (right). In this way, a $3\times 3$ system has to be solved for each pressure grid-point. In particular, the system to solve is as follows,
$$
\begin{pmatrix}
1+\frac{\beta}{d} |u^{\gamma}_{j+1/2}| & 0 & -d\frac{K_f^{\boldsymbol{\tau}}}{h} \\[0.5ex]
0 & 1+\frac{\beta}{d} |u^{\gamma}_{j-1/2}| & d\frac{K_f^{\boldsymbol{\tau}}}{h} \\[0.5ex]
\frac{1}{h} & -\frac{1}{h} & 0
\end{pmatrix}
\begin{pmatrix}
\delta u^{\gamma}_{j+1/2} \\[0.5ex]
\delta u^{\gamma}_{j-1/2} \\[0.5ex]
\delta p^{\gamma}_{j}
\end{pmatrix} = 
\begin{pmatrix}
\delta r^{u^{\gamma}}_{j+1/2} \\[0.5ex]
\delta r^{u^{\gamma}}_{j-1/2} \\[0.5ex]
\delta r^{p^{\gamma}}_{j}
\end{pmatrix},
$$
where the diagonal elements of the matrix are computed by using the last updated values of the velocities, $h$ again denotes the uniform mesh-size, and $K_f^{\boldsymbol{\tau}}$ is assumed constant along the fracture. 
Notice that no coupling between the porous matrix and the fracture unknowns is needed in the smoothing procedure. This is different to the case in which Darcy's law was considered both in the bulk and in the fracture \cite{arr:gas:por:rod:18}. In fact, in order to obtain a robust monolithic multigrid solver in the case of considering barrier type fractures modeled by Darcy's law, we needed to couple both types of unknowns.

\section{Numerical experiments}\label{sec:experiments}

\begin{figure}[t!]
	\begin{center}
		\begin{tabular}{cc}
			\hspace{-0.6cm}\includegraphics[scale = 0.43]{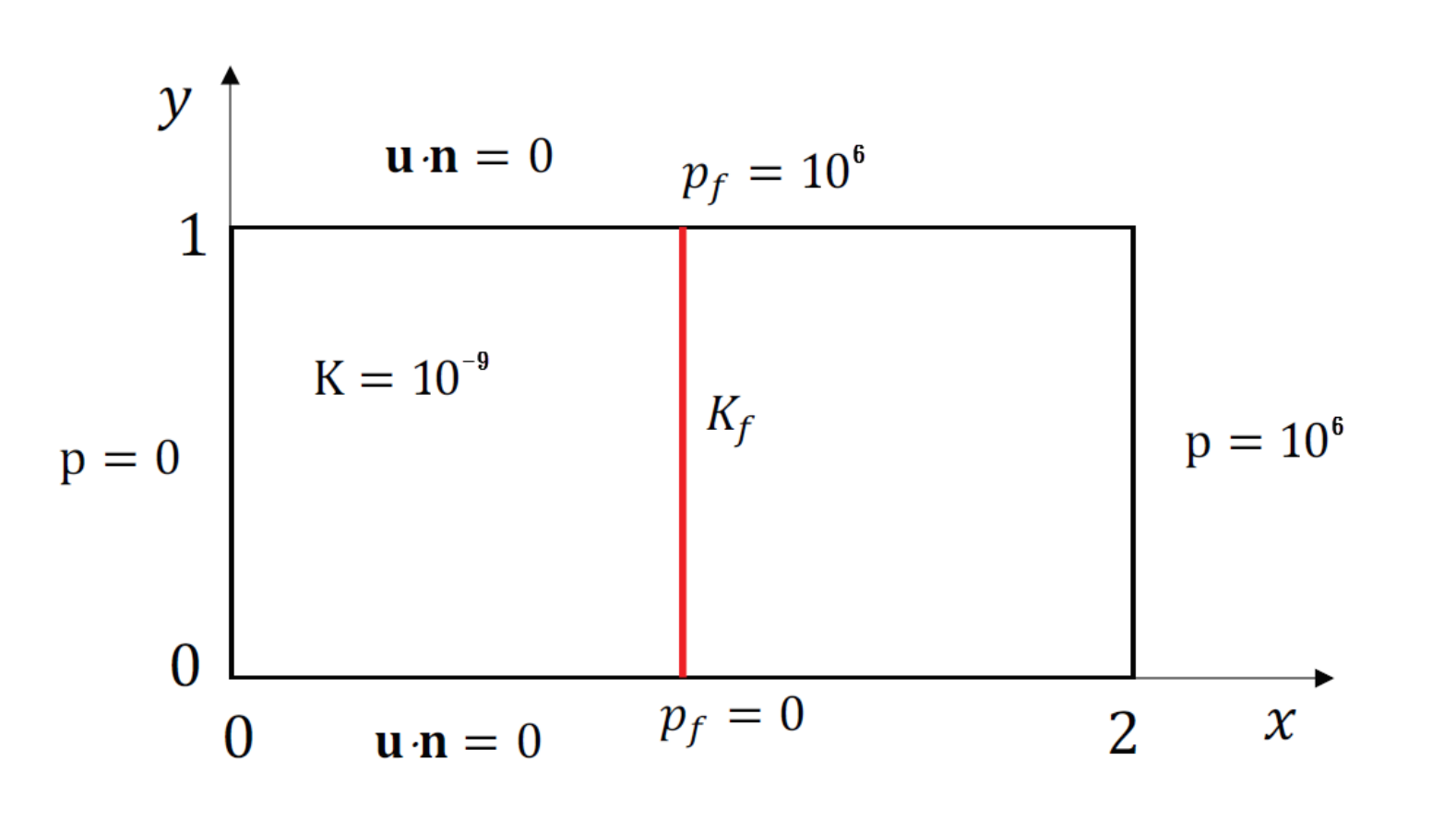}
			&
			\hspace{-0.6cm}\includegraphics[scale = 0.4]{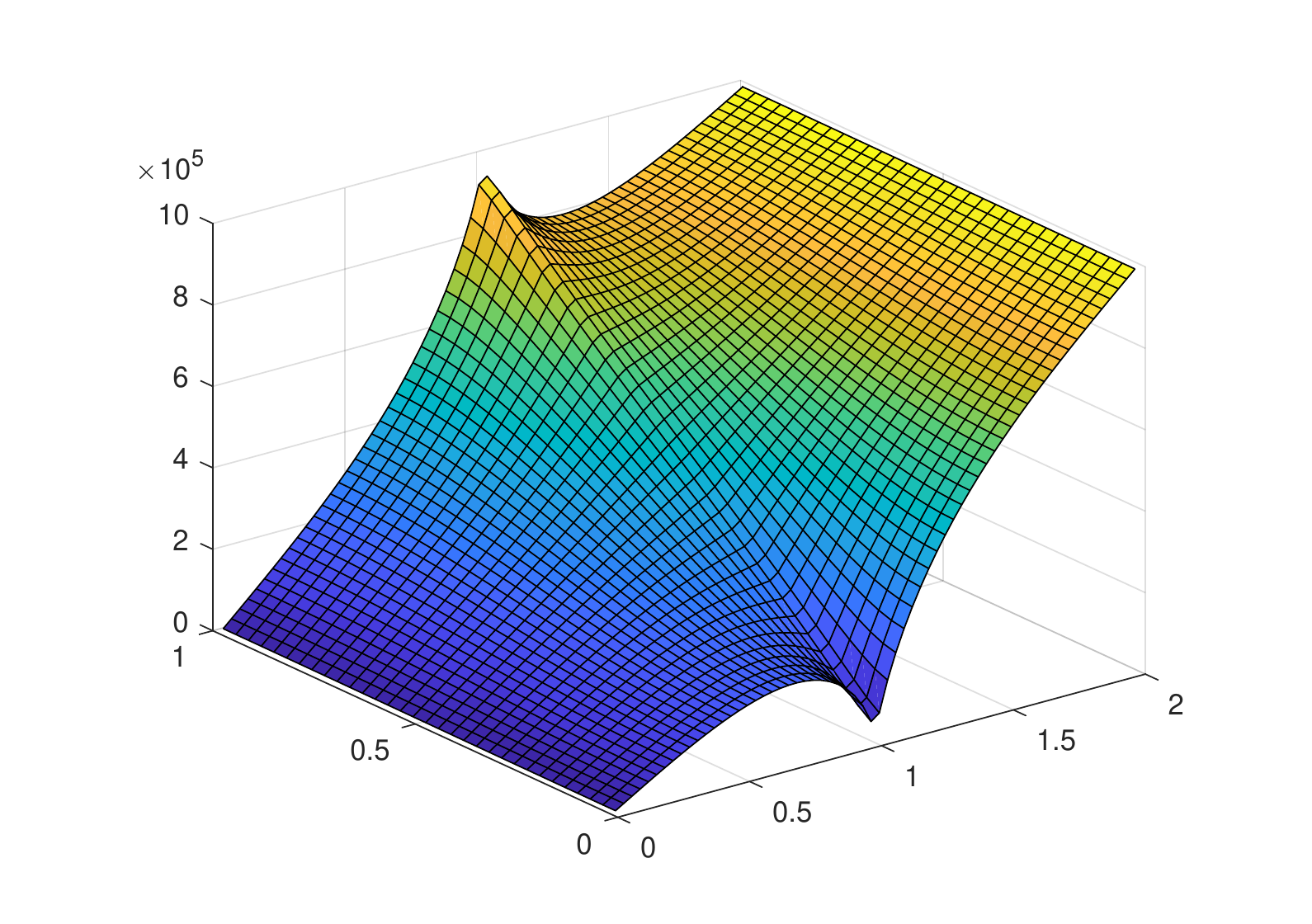}\\
			\hspace{-0.6cm} (a) & \hspace{-0.6cm} (b)
		\end{tabular}
	\end{center}\vspace*{-0.4cm}
	\caption{(a) Domain and boundary conditions and (b) pressure solution for the case $K_f = 10^{-6}$ and $b=10$.}
	\label{first_numerical}
\end{figure}

In this section we present some numerical experiments to show the robustness and the efficiency of the proposed multigrid method for the interface model. We consider a test 
problem presented in \cite{fri:rob:saa:08} where the domain consists of an horizontal 
rectangular slice of porous medium $\Omega = (0,2) \times (0,1)$. Such a domain is divided
into two equally sized subdomains by a vertical fracture $\Omega_f$ of unit length and of width 
$d = 0.01$. The permeability in the porous medium is assumed to be $\mathbf{K}=K\mathbf{I}$, where $K = 10^{-9}$ and $\mathbf{I}$ stands for the identity matrix. In turn, the permeability in the fracture is given by $\mathbf{K}_f=K_f\mathbf{I}$, where $K_f={K}_{f}^{\boldsymbol{\tau}}={K}_{f}^{\mathbf{n}}$ is supposed to be greater than $K$. We will perform several numerical tests for different values of permeability $K_f$ as well as for different values of the Forchheimer number $\beta$. The upper
and lower boundaries of the porous medium are assumed to be impermeable. Pressure is fixed 
on the left boundary to $p=0$ whereas on the right boundary the pressure is fixed to $p=10^6$.
The boundary conditions of the fracture are Dirichlet. More concretely, $p_f = 10^6$ on the top
extremity of the fracture and $p_f=0$ on the bottom. In Figure \ref{first_numerical} (a) we display all these settings. As an example, in Figure \ref{first_numerical} (b) we show the pressure solution obtained when a permeability of $K_f = 10^{-6}$ is assumed in the fracture  and a Forchheimer coefficient 
$\beta= 10$ is considered. The problem has been discretized by the finite volume scheme described in Section \ref{sec:discretization} by considering a uniform grid in both directions with mesh-size $h$.

In all numerical tests, we apply the FAS multigrid method based on Vanka-type smoother described in Section \ref{sec:smoother}. For each smoothing step we will apply one iteration of the Vanka-type method on the porous medium and two iterations of the same method on the fracture, with a damping parameter $w= 0.7$. Notice that in this model the fracture is considered as a one-dimensional object, and therefore the computational cost of solving the fracture unknowns is negligible compared to that corresponding to solving the problem on the porous medium. We use W-cycles with two pre- and two post-smoothing steps. We have seen that this choice gives very good results for solving difficult coupled problems like the Darcy--Stokes system \cite{PeiyaoSISC}, the Biot--Stokes system \cite{PeiyaoJCP} and the single phase Darcy--Darcy coupling between the fractures and the
porous matrix \cite{arr:gas:por:rod:18}. All our numerical computations were carried out using MATLAB.

Throughout this section, we show the robustness of the monolithic mixed-dimensional multigrid method with respect to the spatial discretization parameter $h$, with respect to the permeability of the fracture $K_f$ and the Forchheimer coefficient $\beta$. Notice that the problem 
becomes harder to solve as the permeability of the fracture $K_f$ increases, mainly because
of the big jump between the permeabilities of the fracture and the porous medium. Moreover,
the solution of this multi-dimensional coupled problem becomes more difficult as the coefficient $\beta$ increases. Notice that the Forchheimer coefficient enhances the nonlinearity.  

\begin{figure}[t]
	\begin{center}
		\includegraphics[width = 0.5\textwidth]{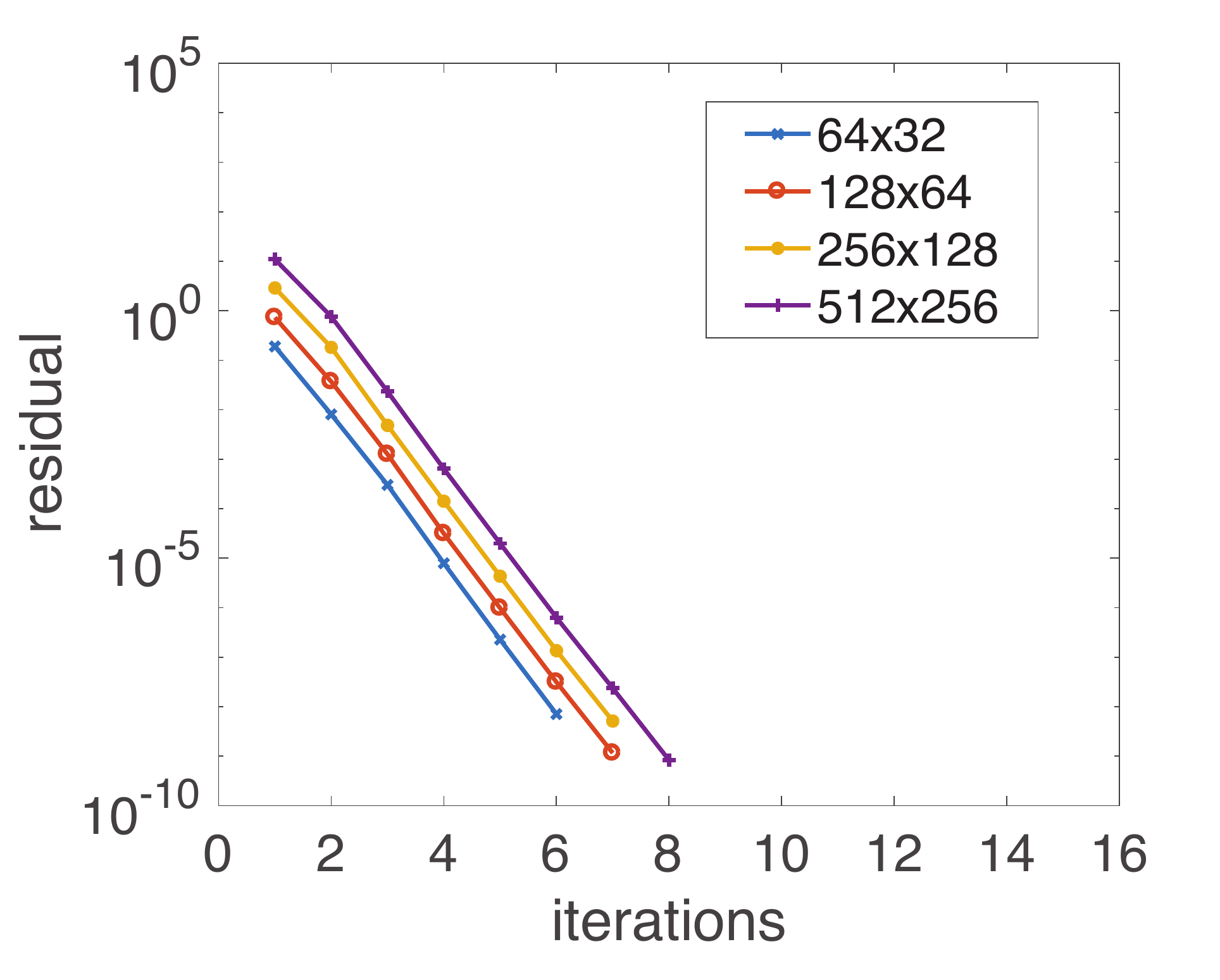}
	\end{center}\vspace*{-0.4cm}
	\caption{History of the convergence of the monolithic multigrid method for $K_f = 10^{-6}$ and $\beta= 10$.}
	\label{figure_convergence}
\end{figure}

We study the performance of the mixed-dimensional multigrid method by fixing the 
permeability of the fracture $K_f = 10^{-6}$ and the Forchheimer coefficient $\beta= 10$. In 
Figure \ref{figure_convergence} we display the history of the convergence of the multigrid solver for different mesh sizes. More concretely, the reduction of the residual is depicted against the number of 
iterations, and the stopping criterium is to reduce the initial residual until $10^{-8}$. It can be 
observed that the convergence of the monolithic mixed-dimensional multigrid method is
independent on the spatial discretization parameter. Moreover, it results in a very efficient 
solver since only around eight iterations are enough to solve this non-linear multidimensional
coupled problem.

Next, we fix the Forchheimer coefficient $\beta= 10$ in order to study the robustness of the mixed-dimensional multigrid method with respect to different values of the permeability of the fracture $K_f$. In Table \ref{table_nonlinear_fixed_kf} we display the number of iterations needed to reduce the initial residual in a factor of $10^{-10}$ for different grid sizes and for different
permeabilities. We can observe that for all values of $K_f$ the performance of the
multigrid method is independent on the spatial discretization parameter, being also
necessary few iterations to reach the stopping criterium.

\begin{table}[t]
	\begin{center}\small
		\begin{tabular}{c||cccc}
			\hline\\[-0.3cm]
			$K_f$ & $h^{-1} = 32$ & $h^{-1} = 64$ & $h^{-1} = 128$ & $h^{-1} = 256$ \\
			\hline\\[-0.3cm]
			$10^{-6}$ & 8   & 8   & 8   & 9 \\
			$10^{-4}$ & 9   & 9   & 9   & 9 \\
			$10^{-2}$ & 9   & 9   & 9   & 10 \\
			$1$          & 10 & 10 & 10 & 11\\
			\hline
		\end{tabular}\vspace*{0.2cm}
		\caption{Number of $W(2,2)$-iterations of the FAS multigrid method required to reduce the initial residual in a factor of $10^{-10}$ for different values of the permeability in the fracture $K_f$ and for different grid-sizes. The Forchheimer coefficient is $\beta= 10$.}
		\label{table_nonlinear_fixed_kf}
	\end{center}
\end{table}

Finally, in Table \ref{table_nonlinear_fixed_b} we show the number of iterations required to
reduce the initial residual in a factor of $10^{-10}$ for different values of the Forchheimer coefficient $\beta$ and different grid sizes. Here, we have fixed the permeability of the fracture as $K_f = 10^{-6}$. Notice that parameter $\beta$ controls the strength of the nonlinearity and the bigger $\beta$, the harder the problem becomes. For comparison, we have also included the case $\beta=0$ which corresponds to consider the Darcy's law in the fracture. It is well known that the FAS scheme for linear problems is theoretically equivalent to the usual linear multigrid scheme
\cite{TOS01}. In this case, the resulting multigrid method is similar to that proposed 
in \cite{arr:gas:por:rod:18} for solving the single-phase Darcy flow in a fractured porous medium. We observe that the performance of the solver is very similar, demonstrating that the proposed multigrid method is also robust with respect to the Forchheimer coefficient $\beta$. 

\begin{table}[t]
	\begin{center}\small
		\begin{tabular}{c||cccc}
			\hline\\[-0.3cm]
			$\beta$ & $h^{-1} = 32$ & $h^{-1} = 64$ & $h^{-1} = 128$ & $h^{-1} = 256$ \\
			\hline
			0     & 8 & 8 & 8   &  8 \\
			10   & 8 & 8 & 8   &  9 \\
			50   & 8 & 8 & 9   & 10 \\
			100 & 8 & 9 & 10 &  10\\
			200 & 9 & 9 & 10 &  10\\
			\hline
		\end{tabular}\vspace*{0.2cm}
		\caption{Number of $W(2,2)$-iterations of the FAS multigrid method required to reduce the initial residual in a factor of $10^{-10}$ for different values of the Forchheimer coefficient $\beta$ and for different grid-sizes. The permeability in the fracture is $K_f = 10^{-6}$.}
		\label{table_nonlinear_fixed_b}
	\end{center}
\end{table}


%
%




\section*{Acknowledgements}
Francisco J. Gaspar has received funding from the European Union's Horizon 2020 research and innovation programme under the Marie Sk\l{}odowska-Curie grant agreement No. 705402, POROSOS. The work of Andr\'es Arrar\'as, Laura Portero and Carmen Rodrigo is supported in part by the FEDER/MINECO project MTM2016-75139-R.




\bibliographystyle{elsarticle-num}
\bibliography{references}







\end{document}